\documentclass[a4paper,11pt]{amsart}
\usepackage{graphicx}
\usepackage{mathptmx}
\usepackage{amsmath}
\usepackage{amssymb}
\usepackage{enumitem}
\usepackage{xcolor}
\usepackage{pgfplots}

\newmuskip\pFqmuskip

\newcommand*\pFq[6][8]{%
  \begingroup % only local assignments
  \pFqmuskip=#1mu\relax
  \mathcode`=\string"8000
  \begingroup\lccode`\~=`\,
  \lowercase{\endgroup\let~}\pFqcomma
  F^{#2}_{#3}{\left(\genfrac..{0pt}{}{#4}{#5}\bigg|#6\right)}%
  \endgroup
}
\newcommand{\pFqcomma}{\mskip\pFqmuskip}

\newtheorem{theorem}{Theorem}[section]

\newtheorem{corollary}[theorem]{Corollary}

\begin{document}

\title[Study on Morgan-Voyce type polynomials with Euler-Seidel algorithm]{Study on Morgan-Voyce type polynomials with Euler-Seidel algorithm}

\author{Taekyun  Kim}
\address{Department of Mathematics, Kwangwoon University, Seoul 139-701, Republic of Korea}
\email{tkkim@kw.ac.kr}
\author{Dae San  Kim}
\address{Department of Mathematics, Sogang University, Seoul 121-742, Republic of Korea}
\email{dskim@sogang.ac.kr}

\subjclass[2010]{11B73; 11B83}
\keywords{Morgan-Voyce type polynomials; Euler-Seidel matrix; degenerate Euler-Seidel matrix}

\begin{abstract}
This paper bridges the domains of degenerate special polynomials, the Euler–Seidel matrix method, and Morgan–Voyce polynomials. We introduce two new families of Morgan–Voyce type polynomials, $M_{n}(x)$ and $N_{n}(x)$, and establish their structural properties, including explicit formulas and recurrence relations. Additionally, we define three polynomial variants, $K_{n}(x)$, $J_{n}(x)$, and $L_{n}(x)$, prove that three distinct binomial-type sums for Bell polynomials can be expressed as finite sums involving these new families, and derive their exponential generating functions. We then construct Euler–Seidel matrices using the initial sequences given by $M_{n}(x)$, $N_{n}(x)$, $K_{n}(x)$, $J_{n}(x)$, and $L_{n}(x)$. Finally, we introduce the degenerate variants $K_{n,\lambda}(x)$, $J_{n,\lambda}(x)$, and $L_{n,\lambda}(x)$, detailing their structural interrelations, their exponential generating  functions and their connection to the degenerate Euler–Seidel matrix framework. Our results yield novel algebraic identities and expand the application of matrix methods in combinatorial analysis.
\end{abstract}

\maketitle

\section{Introduction}
In recent years, the study of degenerate versions of special polynomials and numbers has attracted significant attention from many mathematicians. Since Carlitz introduced the degenerate Bernoulli and Eulerian polynomials \cite{4, 5}, various degenerate variants of special polynomials and numbers have been actively investigated. These variants not only possess interesting combinatorial properties but also provide deep insights into analytic number theory, mathematical physics, and combinatorics \cite{7, 9, 16, 17}. In particular, degenerate special numbers and polynomials--such as the degenerate exponential functions, degenerate logarithms \cite{8}, and degenerate Stirling numbers \cite{14}--are frequently explored using tools like generating functions, $p$-adic analysis, and combinatorial identities \cite{3, 6, 19, 23}. Recent advancements have also incorporated structural and operator-based methodologies, such as Boson operators and Spivey-type recurrence relations, to further characterize degenerate Bell and related polynomials \cite{10, 15}. \par
Concurrently, the Euler-Seidel matrix method has proven to be a remarkably powerful and elegant tool for investigating properties of sequences, generating functions, and polynomials. By arranging an initial sequence into the first row or column of a matrix and applying a specific recurrence relation, one can determine a final sequence. This matrix approach allows researchers to discover underlying algebraic structures and derive explicit formulas, recurrence relations, and convolution identities that might otherwise remain obscured. Notably, this framework has recently been extended to the degenerate domain, establishing algorithmic approaches and matrix methods specifically tailored for degenerate Bernoulli, Euler, and Genocchi polynomials \cite{11, 12, 13}. \par
Another fascinating area of study involves Morgan-Voyce polynomials, which originally arose from the analysis of electrical ladder networks \cite{1, 2, 18}. Over time, they have gained prominence in combinatorics and number theory due to their close relationships with Fibonacci numbers, Lucas numbers, and various orthogonal polynomials. Subsequent generalizations and investigations into their rising diagonals have highlighted their deep algebraic utility \cite{20, 21, 22}. Exploring the intersections between these classical polynomials, their modern degenerate counterparts, and the transforming power of the Euler-Seidel matrix opens up new avenues for discovering deep algebraic identities. \par
Motivated by these ongoing developments, the primary objective of this paper is to bridge these domains. We introduce new families of Morgan-Voyce type polynomials, investigate their structural properties, and apply both standard and degenerate Euler-Seidel matrices to derive novel binomial-type sums involving Bell polynomials and (degenerate) special polynomials. \par
The remainder of this paper is organized as follows. \par
$\indent\bullet$ In Section 2, we review the preliminary concepts necessary for our discussion. Specifically, we recall the degenerate exponentials, degenerate logarithms, unsigned degenerate Stirling numbers of the first kind, degenerate Stirling numbers of the second kind, and the Bell polynomials. We also remind the reader of the standard Euler-Seidel matrix and the degenerate Euler-Seidel matrix, outlining the relations between their initial and final sequences. Equivalently, we state these relations in terms of the exponential generating functions of their respective initial and final sequences. Finally, we introduce the two classic Morgan-Voyce polynomials. \par
$\indent\bullet$ Section 3 presents the main results of this paper, beginning with the introduction of two new Morgan-Voyce type polynomials, denoted by $M_{n}(x)$ and $N_{n}(x)$. In Theorem 3.1, we establish a relation between $N_{n}(x)$ and $M_{n}(x)$, derive a recurrence relation for $M_{n}(x)$, and provide an explicit expression for $M_{n}(x)$. An explicit expression for $N_{n}(x)$ is subsequently provided in Corollary 3.3. We then explore the Euler-Seidel matrices associated with the initial sequences given by these two Morgan-Voyce type polynomials. This leads us to define two new types of Morgan-Voyce polynomials: $K_{n}(x)$ corresponding to $M_{n}(x)$, and $J_{n}(x)$ corresponding to $N_{n}(x)$. Additionally, we introduce another polynomial, $L_{n}(x)$, derived from $J_{n}(x)$. In Theorems 3.4, 3.6, and 3.7, we demonstrate that three distinct binomial-type sums for Bell polynomials are respectively equal to finite sums involving $K_{n}(x)$, $J_{n}(x)$, and $L_{n}(x)$. Connections between these new polynomials--specifically the relations between $K_{n}(x)$ and $J_{n}(x)$, and between $J_{n}(x)$ and $L_{n}(x)$--are determined in Theorems 3.5 and 3.8, respectively. Concluding this section, we examine the Euler-Seidel matrices associated with the initial sequences given by $K_{n}(x)$, $J_{n}(x)$, and $L_{n}(x)$. \par
$\indent\bullet$ In Section 4, we introduce the polynomials $K_{n,\lambda}(x)$, $J_{n,\lambda}(x)$, and $L_{n,\lambda}(x)$, which serve as the degenerate versions of $K_{n}(x)$, $J_{n}(x)$, and $L_{n}(x)$, respectively. We establish a relation between $K_{n,\lambda}(x)$ and $J_{n,\lambda}(x)$ in Theorem 4.1, as well as a relation between $J_{n,\lambda}(x)$ and $L_{n,\lambda}(x)$ in Theorem 4.2. Finally, we analyze the Euler-Seidel matrix associated with the initial sequence given by $L_{n,\lambda}(x)$.

\section{Preliminaries}
Here we review the preliminary concepts necessary for our discussion. \\
For any nonzero $\lambda\in\mathbb{R}$, the degenerate exponentials are defined by \cite{8,9}
\begin{equation}
e_{\lambda}^{x}(t)=\sum_{k=0}^{\infty}(x)_{k,\lambda}\frac{t^{k}}{k!},\quad e_{\lambda}(t)=e_{\lambda}^{1}(t), \label{1}
\end{equation}
where
\begin{equation*}
(x)_{0,\lambda}=1,\quad (x)_{n,\lambda}=x(x-\lambda)(x-2\lambda)\cdots\big(x-(n-1)\lambda\big),\ (n\ge 1).
\end{equation*}
Note that
\begin{equation*}
\lim_{\lambda\rightarrow 0}e_{\lambda}^{x}(t)=e^{xt}=\sum_{n=0}^{\infty}\frac{x^{n}}{n!}t^{n}.
\end{equation*}
Let $\log_{\lambda}t$ be the compositional inverse of $e_{\lambda}(t)$, called the degenerate logarithms, such that \cite{8,9}
\begin{equation*}
\log_{\lambda}\big(e_{\lambda}(t)\big)=e_{\lambda}\big(\log_{\lambda}(t)\big)=t.
\end{equation*}
Then we have
\begin{equation}
\log_{\lambda}(1+t)=\sum_{k=0}^{\infty}\frac{\lambda^{k-1}(1)_{k,1/\lambda}}{k!}t^{k}=\frac{1}{\lambda}\big((1+t)^{\lambda}-1\big).\label{2}
\end{equation}
The unsigned degenerate Stirling numbers of the first kind are defined by \cite{8,9}
\begin{equation}
\frac{1}{k!}\big(-\log_{\lambda}(1-t)\big)^{k}=\sum_{n=k}^{\infty}{n \brack k}_{\lambda}\frac{t^{n}}{n!},\ (k\ge 0). \label{3}
\end{equation}
The degenerate Stirling numbers of the second kind are given by \cite{8,9}
\begin{equation}
\frac{1}{k!}\big(e_{\lambda}(t)-1\big)^{k}=\sum_{n=k}^{\infty}{n \brace k}_{\lambda}\frac{t^{n}}{n!},\quad (k\ge 0).\label{4}
\end{equation}
Note that
\begin{align*}
&\lim_{\lambda\rightarrow 0}\frac{1}{k!}\big(-\log_{\lambda}(1-t)\big)^{k}=\frac{1}{k!}\log^{k}\bigg(\frac{1}{1-t}\bigg)=\sum_{n=k}^{\infty}{n \brack k}\frac{t^{n}}{n!}, \\
&\lim_{\lambda\rightarrow 0}\frac{1}{k!}\big(e_{\lambda}(t)-1\big)^{k}=\frac{1}{k!}\big(e^{t}-1\big)^{k}=\sum_{n=k}^{\infty}{n \brace k}\frac{t^{n}}{n!},
\end{align*}
where ${n \brack k}$ and ${n \brace k}$ are the unsigned Stirling numbers of the first kind and the Stirling numbers of the second kind. \par
It is well known that the Bell polynomials are defined by \cite{6,19}
\begin{equation}
\phi_{n}(x)=\sum_{k=0}^{n}{n \brace k}x^{k},\quad (n\ge 0). \label{5}
\end{equation}
Thus, by \eqref{5}, we easily get
\begin{equation}
e^{x(e^{t}-1)}=\sum_{n=0}^{\infty}\phi_{n}(x)\frac{t^{n}}{n!}.\label{6}
\end{equation} \par
It is known that the Euler-Seidel matrix $(a_{n,k}(x))_{n,k\ge 0}$ is defined by \cite{11,12}
\begin{equation}
\begin{aligned}
&a_{0,n}(x)=a_{n}(x),\quad (n\ge 0), \\
&a_{k,n}(x)=a_{k-1,n}(x)+a_{k-1,n+1}(x),\quad (k\ge 1,n\ge 0).
\end{aligned}\label{7}
\end{equation}
The Euler-Seidel matrix associated with $(a_{n,k}(x))_{n,k\ge 0}$ is given by
\begin{equation*}
A=\big(a_{n,k}(x)\big)_{n,k\ge 0}=\left(\begin{array}{ccccc}
a_{0,0}(x) & a_{0,1}(x) & a_{0,2}(x) & a_{0,3}(x) & \cdots \\
a_{1,0}(x) & a_{1,1}(x) & a_{1,2}(x) & a_{1,3}(x) & \cdots \\
a_{2,0}(x) & a_{2,1}(x) & a_{2,2}(x) & a_{2,3}(x) & \cdots \\
\vdots & \vdots & \vdots & \vdots & \cdots
\end{array}\right).
\end{equation*}
Thus, by \eqref{7}, we easily get
\begin{equation}
a_{n,0}(x)=\sum_{k=0}^{n}\binom{n}{k}a_{0,k}(x)\ \Longleftrightarrow\ a_{0,n}(x)=\sum_{k=0}^{n}\binom{n}{k}(-1)^{n-k}a_{k,0}(x). \label{8}
\end{equation} \par
In \cite{11,12}, the degenerate Euler-Seidel matrix $(a_{n,k}(x|\lambda))_{n,k\ge 0}$ is given by
\begin{equation}
\begin{aligned}
&a_{0,n}(x|\lambda)=a_{n}(x|\lambda),\\
&a_{k,n}(x|\lambda)=\big(1-(k-n)\lambda\big)a_{k-1,n}(x|\lambda)+a_{k-1,n+1}(x|\lambda),
\end{aligned} \label{9}
\end{equation}
where $k\ge 1$ and $n\ge 0$. \par
Note that
\begin{equation}
a_{n,0}(x|\lambda)=\sum_{k=0}^{n}\binom{n}{k}(1-\lambda)_{n-k,\lambda}a_{0,k}(x|\lambda),\quad (n\ge 0). \label{10}
\end{equation}
The degenerate Euler-Seidel matrix associated with $\big(a_{0,n}(x|\lambda)\big)_{n\ge 0}$ is given by
\begin{equation*}
A=\big(a_{n,k}(x|\lambda)\big)_{n,k\ge 0}=\left(\begin{matrix}
a_{0,0}(x|\lambda) & a_{0,1}(x|\lambda) & a_{0,2}(x|\lambda) & \cdots \\
a_{1,0}(x|\lambda) & a_{1,1}(x|\lambda) & a_{1,2}(x|\lambda) & \cdots \\
a_{2,0}(x|\lambda) & a_{2,1}(x|\lambda) & a_{2,2}(x|\lambda) & \cdots\\
\vdots & \vdots & \vdots & \cdots
\end{matrix}\right).
\end{equation*} \par
Let $\displaystyle A(t)=\sum_{n=0}^{\infty}a_{0,n}(x)\frac{t^{n}}{n!}\displaystyle$. Then \eqref{8} is equivalent to \cite{11,12}
\begin{equation}
\overline{A}(t)=\sum_{n=0}^{\infty}a_{n,0}(x)\frac{t^{n}}{n!}=e^{t}A(t)\quad \mathrm{(Seidel)}.\label{11}
\end{equation}
In addition, if
\begin{equation*}
A_{\lambda}(t)=\sum_{n=0}^{\infty}a_{0,n}(x|\lambda)\frac{t^{n}}{n!}\quad \mathrm{(degenerate\,\, Seidel)},
\end{equation*}
then \eqref{10} is equivalent to \cite{11,12}
\begin{equation}
\overline{A}_{\lambda}(t)=\sum_{n=0}^{\infty}a_{n,0}(x|\lambda)\frac{t^{n}}{n!}=e_{\lambda}^{1-\lambda}(t)A_{\lambda}(t). \label{12}
\end{equation} \par
Two polynomials $R_{n}(x)$ and $b_{n}(x)$ of Morgan-Voyce \cite{1,2,18} are defined by
\begin{equation}
b_{n}(x)=xB_{n-1}(x)+b_{n-1}(x),\quad B_{n}(x)=(x+1)B_{n-1}(x)+b_{n-1}(x),\quad (n\ge 1),\label{13}
\end{equation}
where $b_{0}=B_{0}(x)=1$. \par

\section{Study on Morgan-Voyce type polynomials with Euler-Seidel algorithm}
In this section, we introduce and several new Morgan–Voyce type polynomials, namely $M_{n}(x)$, $N_{n}(x)$, $K_{n}(x)$, $J_{n}(x)$, and $L_{n}(x)$, deriving their explicit expressions, exponential generating functions and binomial-type sums for Bell polynomials as finite sums involving these polynomials. Furthermore, we construct Euler–Seidel matrices using the initial sequences given by $M_{n}(x)$, $N_{n}(x)$, $K_{n}(x)$, $J_{n}(x)$, and $L_{n}(x)$. \par
Now, we consider two Morgan-Voyce type polynomials given by (see \eqref{13})
\begin{equation}
\begin{aligned}
&N_{n}(x)=nxM_{n-1}(x)+nN_{n-1}(x),\\
&M_{n}(x)=n(x+1)M_{n-1}(x)+nN_{n-1}(x),\quad (n\ge 1),
\end{aligned}\label{14}	
\end{equation}
where $N_{0}(x)=M_{0}(x)=1$. \par
\vspace{0.1in}
From the recurrence relations in \eqref{14}, we find
\begin{align*}
&M_{1}(x)=x+2,\quad M_{2}(x)=2x^{2}+8x+6=2!(x^{2}+4x+3),\\
&M_{3}(x)=6x^{3}+36x^{2}+60x+24=3!(x^{3}+6x^{2}+10x+4),\\
&M_{4}(x)=24x^{4}+192x^{3}+504x^{2}+480x+120=4!(x^{4}+8x^{3}+21x^{2}+20x+5),\\
&N_{1}(x)=x+1,\quad N_{2}(x)=2x^{2}+6x+2=2!(x^{2}+3x+1),\\
&N_{3}(x)=6x^{3}+30x^{2}+36x+6=3!(x^{3}+5x^{2}+6x+1),\\
&N_{4}(x)=24x^{4}+168x^{3}+360x^{2}+240x+24=4!(x^{4}+7x^{3}+15x^{2}+10x+1).
\end{align*}
\vspace{0.1in}
From \eqref{14}, we have
\begin{equation}
N_{n}(x)-nN_{n-1}(x)=nxM_{n-1}(x). \label{15}
\end{equation}
Thus, by \eqref{14} and \eqref{15}, we get
\begin{align}
M_{n}(x)&=nxM_{n-1}(x)+nM_{n-1}(x)+nN_{n-1}(x) \label{16}\\
&=N_{n}(x)-nN_{n-1}(x)+nM_{n-1}(x)+nN_{n-1}(x).\nonumber
\end{align}
From \eqref{16}, we have
\begin{equation}
N_{n}(x)=M_{n}(x)-nM_{n-1}(x),\quad (n\ge 1).\label{17}
\end{equation}
On the other hand, by \eqref{14}, we get
\begin{align}
M_{n+1}(x)&=(n+1)(x+1)M_{n}(x)+(n+1)N_{n}(x)\label{18}\\
&=(n+1)(x+1)M_{n}(x)+(n+1)\big(M_{n}(x)-nM_{n-1}(x)\big)\nonumber \\
&=(n+1)(x+2)M_{n}(x)-n(n+1)M_{n-1}(x).\nonumber
\end{align} \par
Now, we observe that
\begin{align*}
&M_{0}(x)=1,\quad M_{1}(x)=x+2=\sum_{k=0}^{1}\binom{2+k}{1-k}x^{k},\\
&M_{2}(x)=2!(x^{2}+4x+3)=2!\sum_{k=0}^{2}\binom{3+k}{2-k}x^{k}. \nonumber
\end{align*}
Assume that
\begin{equation}
M_{n}(x)=n!\sum_{k=0}^{n}\binom{n+k+1}{n-k}x^{k},\quad (n\ge 1). \label{19}	
\end{equation}
From \eqref{18} and \eqref{19}, we note that
\begin{align*}
M_{n+1}(x)&=(n+1)(x+2)M_{n}(x)-n(n+1)M_{n-1}(x) \\
&=(n+1)xM_{n}(x)+(n+1)M_{n}(x)+(n+1)M_{n}(x)-n(n+1)M_{n-1}(x) \\
&=(n+1)!\sum_{k=0}^{n}\binom{n+k+1}{n-k}x^{k+1}+(n+1)!\sum_{k=0}^{n}\binom{n+k+1}{n-k}x^{k}\\
&\quad+(n+1)!\sum_{k=0}^{n}\binom{n+k+1}{n-k}x^{k}-(n+1)!\sum_{k=0}^{n-1}\binom{n+k}{n-k-1}x^{k}\\
&=(n+1)!\sum_{k=0}^{n+1}\bigg(\binom{n+k}{n-k+1}+\binom{n+k+1}{n-k}\bigg)x^{k}\\
&\quad+(n+1)!\sum_{k=0}^{n}\bigg(\binom{n+k+1}{n-k}-\binom{n+k}{n-k-1}\bigg)x^{k}\\
&=(n+1)!\sum_{k=0}^{n+1}\bigg(\binom{n+k}{n-k+1}+\binom{n+k+1}{n-k}\bigg)x^{k}\\
&\quad+(n+1)!\sum_{k=0}^{n}\binom{n+k}{n-k}x^{k}\\
&=(n+1)!\sum_{k=0}^{n+1}\bigg(\binom{n+k+1}{n-k+1}+\binom{n+k+1}{n-k}\bigg)x^{k}\\
&=(n+1)!\sum_{k=0}^{n+1}\binom{n+2+k}{n+1-k}x^{k},
\end{align*}
which shows that \eqref{19} is valid also for $n+1$.
Therefore, by \eqref{17}, \eqref{18} and \eqref{19}, we obtain the following theorem.
\begin{theorem}
For $n\ge 1$, we have
\begin{equation*}
N_{n}(x)=M_{n}(x)-nN_{n-1}(x),
\end{equation*}
and
\begin{equation*}
M_{n+1}(x)=(n+1)(x+2)M_{n}(x)-n(n+1)M_{n-1}(x).
\end{equation*}
Moreover, we have the following explicit expression
\begin{equation}
M_{n}(x)=n!\sum_{k=0}^{n}\binom{n+k+1}{n-k}x^{k},\quad (n \ge 0).\label{20}
\end{equation}
\end{theorem}
From \eqref{20}, we obtain the next corollary.
\begin{corollary}
For $n,m\ge 1$, we have
\begin{equation*}
M_{n+m}(x)=\binom{n+m}{n}\big(M_{n}(x)M_{m}(x)-nmM_{n-1}(x)M_{m-1}(x)\big).
\end{equation*}
\end{corollary}
From \eqref{17}, we have
\begin{align}
N_{n}(x)&=M_{n}(x)-nM_{n-1}(x)\label{21}\\
&= n!\sum_{k=0}^{n}\binom{n+k+1}{n-k}x^{k}-n(n-1)!\sum_{k=0}^{n-1}\binom{n+k}{n-k-1}x^{k}\nonumber\\
&=n!\sum_{k=0}^{n}\bigg(\binom{n+k+1}{n-k}-\binom{n+k}{n-k-1}\bigg)x^{k} \nonumber\\
&=n!\sum_{k=0}^{n}\binom{n+k}{n-k}x^{k}.\nonumber
\end{align}
From \eqref{21}, we obtain the following corollary.
\begin{corollary}
For $n\ge 0$, we have
\begin{equation*}
N_{n}(x)=n!\sum_{k=0}^{n}\binom{n+k}{n-k}x^{k}.
\end{equation*}
\end{corollary}
First, we let $a_{0,n}(x)=M_{n}(x),\ (n\ge 0)$, in \eqref{2}. Then, by the Seidel's formula in \eqref{12}, we get
\begin{equation*}
a_{n,0}(x)=\sum_{k=0}^{n}\binom{n}{k}M_{k}(x),\quad (n\ge 0). 	
\end{equation*}
The Euler-Seidel matrix associated with $\big(a_{0,n}(x)\big)_{n\ge 0}=\big(M_{n}(x)\big)_{n\ge 0}$ is given by
\begin{align*}
\left(\begin{matrix}
1 & x+2 & 2x^{2}+8x+6 & \cdots \\
x+3 & 2x^{2}+9x+8 & 6x^{3}+38x^{2}+68x+30 & \cdots\\
2x^{2}+10x+11 & 6x^{3}+40x^{2}+77x+38 & \vdots & \cdots \\
6x^{3}+42x^{2}+87x+49 & \vdots & \vdots & \cdots \\
\vdots & \vdots & \vdots & \cdots
\end{matrix}\right).
\end{align*}
Next, we let $a_{0,n}(x)=N_{n}(x),\ (n\ge 0)$, in \eqref{2}. Then, by the Seidel's formula in \eqref{12}, we have
\begin{equation*}
a_{n,0}(x)=\sum_{k=0}^{n}\binom{n}{k}N_{k}(x),\ (n\ge 0).
\end{equation*}
The Euler-Seidel matrix associated with $\big(a_{0,n}(x)\big)_{n\ge 0}=\big(N_{n}(x)\big)_{n\ge 0}$ is given by
\begin{equation*}
\begin{aligned}
\left(\begin{matrix}
1 & x+1 & 2x^{2}+6x+2 & \cdots \\
x+2 & 2x^{3}+7x+3 & 6x^{3}+32x^{2}+42x+8 & \cdots \\
2x^{3}+8x+5 & 6x^{3}+34x^{2}+49x+11  & \vdots & \cdots \\
6x^{3}+36x^{2}+57x+16  & \vdots & \vdots & \cdots \\
\vdots & \vdots & \vdots & \cdots
\end{matrix}\right).
\end{aligned}
\end{equation*}
In view of \eqref{20}, we define the new type Morgan-Voyce polynomials by
\begin{equation}
K_{n}(x)=n!\sum_{k=0}^{n}\binom{n+k+1}{n-k}\frac{x^{k}}{k!},\quad (n\ge 0), \label{22}	
\end{equation}
and
\begin{equation}
J_{n}(x)=K_{n}(x)-nK_{n-1}(x)=n!\sum_{k=0}^{n}\binom{n+k}{n-k}\frac{x^{k}}{k!},\quad (n\ge 1),\label{23}
\end{equation}
where $K_{0}(x)=J_{0}(x)=1$. \par
From \eqref{22}, we note that
\begin{align}
\sum_{n=0}^{\infty}K_{n}(x)\frac{t^{n}}{n!}&=\sum_{n=0}^{\infty}\sum_{k=0}^{n}\binom{n+k+1}{n-k}\frac{x^{k}}{k!}t^{n} \label{24}	\\
&=\sum_{k=0}^{\infty}\frac{x^{k}}{k!}\sum_{n=k}^{\infty}\binom{n+k+1}{n-k}t^{n}=\sum_{k=0}^{\infty}\frac{x^{k}t^{k}}{k!}\sum_{n=0}^{\infty}\binom{n+2k+1}{n}t^{n}\nonumber\\
&=\frac{1}{(1-t)^{2}}\sum_{k=0}^{\infty}\frac{1}{k!}\bigg(\frac{xt}{(1-t)^{2}}\bigg)^{k}=\frac{1}{(1-t)^{2}}e^{\frac{xt}{(1-t)^{2}}}.\nonumber
\end{align}
Thus, by \eqref{24}, we get
\begin{equation}
\frac{1}{(1-t)^{2}}e^{\frac{xt}{(1-t)^{2}}}=\sum_{n=0}^{\infty}K_{n}(x)\frac{t^{n}}{n!}.\label{25}
\end{equation}
Replacing $t$ by $1-e^{-t}$ in \eqref{25} and recalling \eqref{9}, we have
\begin{align}
\sum_{l=0}^{\infty}K_{l}(x)\frac{1}{l!}\big(1-e^{-t}\big)^{l}&=e^{2t}e^{x(e^{2t}-1)} e^{-x(e^{t}-1)} \label{26}\\
&=\frac{1}{2x}\bigg[\frac{d}{dt}e^{x(e^{2t}-1)}\bigg] e^{-x(e^{t}-1)}\nonumber\\
&=\sum_{l=0}^{\infty}\frac{1}{x}\phi_{l+1}(x)2^{l}\frac{t^{l}}{l!}\sum_{k=0}^{\infty}\phi_{k}(-x)\frac{t^{k}}{k!} \nonumber\\
&=\sum_{n=0}^{\infty}\frac{1}{x}\sum_{l=0}^{n}\binom{n}{l}\phi_{l+1}(x)2^{l}\phi_{n-l}(-x)\frac{t^{n}}{n!}.\nonumber
\end{align}
Using the generating function of the Stirling number of the second kind, we also have (see \eqref{7})
\begin{align}
\sum_{l=0}^{\infty}K_{l}(x)\frac{1}{l!}\big(1-e^{-t}\big)^{l}&=\sum_{l=0}^{\infty}K_{l}(x)\frac{(-1)^{l}}{l!}\big(e^{-t}-1\big)^{l} \label{27} 	\\
&=\sum_{l=0}^{\infty}K_{l}(x)(-1)^{l}\sum_{n=l}^{\infty}{n \brace l}(-1)^{n}\frac{t^{n}}{n!} \nonumber \\
&=\sum_{n=0}^{\infty}\sum_{l=0}^{n}{n \brace l}(-1)^{n-l}K_{l}(x)\frac{t^{n}}{n!}. \nonumber
\end{align}
Therefore, by \eqref{26} and \eqref{27}, we obtain the following theorem.
\begin{theorem}
For $n\ge 0$, we have
\begin{equation*}
x\sum_{l=0}^{n}{n \brace l}K_{l}(x)(-1)^{n-l}=\sum_{l=0}^{n}\binom{n}{l}2^{l}\phi_{l+1}(x)\phi_{n-l}(-x).
\end{equation*}
\end{theorem}
We let $a_{0,n}(x)=K_{n}(x),\ (n\ge 0)$, in \eqref{2}. Then, by the Seidel's formula in \eqref{12} and \eqref{25}, we have
\begin{equation*}
a_{n,0}(x)=\sum_{l=0}^{n}\binom{n}{l}K_{l}(x),\quad \sum_{n=0}^{\infty}a_{n,0}(x)\frac{t^n}{n!}=e^{t}\frac{1}{(1-t)^{2}}e^{\frac{xt}{(1-t^{2}}} .
\end{equation*}
From \eqref{22}, we determine $K_{n}(x)$, for the first few values of $n$:
\begin{align*}
&K_0(x) = 1, \quad K_1(x) = x + 2, \quad K_2(x) = x^2 + 8x + 6, \\
&K_3(x) = x^3 + 18x^2 + 60x + 24, \quad
K_4(x) = x^4 + 32x^3 + 252x^2 + 480x + 120. \\
\end{align*}
The Euler-Seidel matrix associated with $(a_{0,n}(x))_{n\ge 0}=(K_{n}(x))_{n\ge 0}$ is given by
\begin{align*}
\left(\begin{matrix}
1 & x+2 & x^{2}+8x+6 & \cdots \\
x+3 & x^{2}+9x+8 & x^{3}+19x^{2}+68x+30 & \cdots \\
x^{2}+10x+11 & x^{3}+20x^{2}+77x+38 & \vdots & \cdots\\
x^{3}+21x^{2}+87x+49 & \vdots & \vdots & \cdots \\
\vdots & \vdots & \vdots & \cdots
\end{matrix}\right).
\end{align*}
From \eqref{22} and \eqref{23}, we note that
\begin{align}
K_{n}(x)&=n!\sum_{j=0}^{n}\binom{n+j+1}{n-j}\frac{x^{j}}{j!}=n!\sum_{j=0}^{n}\frac{x^{j}}{j!}\binom{n+j+1}{2j+1} \label{28}\\
&=n!\sum_{j=0}^{n}\frac{x^{j}}{j!}\sum_{k=j}^{n}\bigg(\binom{k+j+1}{2j+1}-\binom{k+j}{2j+1}\bigg)=n!\sum_{j=0}^{n}\frac{x^{j}}{j!}\sum_{k=j}^{n}\binom{k+j}{2j}\nonumber\\
&=n!\sum_{j=0}^{n}\frac{x^{j}}{j!}\sum_{k=j}^{n}\binom{k+j}{k-j}=n!\sum_{k=0}^{n}\frac{1}{k!}k!\sum_{j=0}^{k}\binom{k+j}{k-j}\frac{x^{j}}{j!}\nonumber\\
&=n!\sum_{k=0}^{n}\frac{1}{k!}J_{k}(x).\nonumber
\end{align}
Therefore, by \eqref{28}, we obtain the following theorem.
\begin{theorem}
For $n\ge 0$, we have
\begin{equation*}
K_{n}(x)=n!\sum_{k=0}^{n}\frac{J_{k}(x)}{k!}.
\end{equation*}
\end{theorem}
From \eqref{23}, we note that
\begin{align}
\sum_{n=0}^{\infty}\frac{t^{n}}{n!}J_{n}(x)&=\sum_{n=0}^{\infty}t^{n}\sum_{k=0}^{n}\binom{n+k}{n-k}\frac{x^{k}}{k!}=\sum_{k=0}^{\infty}\frac{x^{k}}{k!}\sum_{n=k}^{\infty}\binom{n+k}{n-k}t^{n} \label{29}\\
&=\sum_{k=0}^{\infty}\frac{(tx)^{k}}{k!}\sum_{n=0}^{\infty}\binom{n+2k}{n}t^{n}=\sum_{k=0}^{\infty}\frac{(tx)^{k}}{k!}\bigg(\frac{1}{1-t}\bigg)^{2k+1}\nonumber\\
&=\frac{1}{1-t}e^{\frac{xt}{(1-t)^{2}}}.\nonumber
\end{align}
Thus, by \eqref{29}, we get
\begin{equation}
\frac{1}{1-t}e^{\frac{xt}{(1-t)^{2}}}=\sum_{n=0}^{\infty}J_{n}(x)\frac{t^{n}}{n!}. \label{30}
\end{equation}
Replacing $t$ by $1-e^{-t}$ in \eqref{30}, we have
\begin{align}
\sum_{l=0}^{\infty}J_{l}(x)\frac{1}{l!}\big(1-e^{-t}\big)^{l}&=e^{t}e^{x(e^{2t}-1)}e^{-x(e^{t}-1)}\label{31} \\
&=e^{x(e^{2t}-1)}\bigg[-\frac{1}{x}\frac{d}{dt}e^{-x(e^{t}-1)}\bigg] \nonumber\\
&=\sum_{l=0}^{\infty}\phi_{l}(x)2^{l}\frac{t^{l}}{l!}\bigg(-\frac{1}{x}\sum_{m=0}^{\infty}\phi_{m+1}(-x)\frac{t^{m}}{m!}\bigg)\nonumber\\
&=-\frac{1}{x}\sum_{n=0}^{\infty}\sum_{l=0}^{n}\binom{n}{l}2^{l}\phi_{l}(x)\phi_{n-l+1}(-x)\frac{t^{n}}{n!}. \nonumber
\end{align}
Proceeding in the same way as in \eqref{27}, we also have
\begin{equation}
\sum_{l=0}^{\infty}J_{l}(x)\frac{1}{l!}\big(1-e^{-t}\big)^{l}=\sum_{n=0}^{\infty}\sum_{l=0}^{n}{n \brace l}(-1)^{n-l}J_{l}(x)\frac{t^{n}}{n!}. \label{32}
\end{equation}
Therefore, by \eqref{31} and \eqref{32}, we obtain the following theorem.
\begin{theorem}
For $n\ge 0$, we have
\begin{equation*}
x\sum_{l=0}^{n}{n \brace l}J_{l}(x)(-1)^{n-l-1}=\sum_{l=0}^{n}\binom{n}{l}2^{l}\phi_{l}(x)\phi_{n+1-l}(-x).
\end{equation*}
\end{theorem}
We let $a_{0,n}(x)=J_{n}(x),\ (n\ge 0)$, in \eqref{2}. Then, by the Seidel's formula in \eqref{12} and \eqref{30}, we get
\begin{equation*}
a_{n,0}(x)=\sum_{k=0}^{n}\binom{n}{k}J_{k}(x), \quad \sum_{n=0}^{\infty}a_{n,0}(x)\frac{t^n}{n!}=e^{t}\frac{1}{1-t}e^{\frac{xt}{(1-t^{2}}}. 	
\end{equation*}
From \eqref{23}, we determine $J_{n}(x)$, for the first few values of $n$:
\begin{align*}
&J_0(x) = 1, \quad J_1(x) = x + 1, \quad J_2(x) = x^2 + 6x + 2, \\
&J_3(x) = x^3 + 15x^2 + 36x + 6, \quad
J_4(x) = x^4 + 28x^3 + 180x^2 + 240x + 24. \\
\end{align*}
The Euler-Seidel matrix associated with $\big(a_{0,n}(x)\big)_{n\ge 0}=\big(J_{n}(x)\big)_{n\ge 0}$ is given by
\begin{align*}
\left(\begin{matrix}
1 & x+1 & x^{2}+6x+2 & \cdots \\
x+2 & x^{2}+7x+3 & x^{3}+16x^{2}+42x+8 & \cdots \\
x^{2}+8x+5 & x^{3}+17x^{2}+49x+11& \vdots & \cdots \\
x^{3}+18x^{2}+57x+16 & \vdots & \vdots & \cdots \\
\vdots & \vdots & \vdots & \cdots \\
\end{matrix}\right).
\end{align*}
We consider the new polynomials $L_{n}(x)$, arising from $J_{n}(x)$, which are given by
\begin{equation}
L_{n}(x)=J_{n}(x)-nJ_{n-1}(x),\ (n\ge 0).\label{33}
\end{equation}
Then, by \eqref{23} and \eqref{33}, we get
\begin{align}
L_{n}(x)&=J_{n}(x)-nJ_{n-1}(x)=n!\sum_{k=0}^{n}\binom{n+k}{n-k}\frac{x^{k}}{k!}-n!\sum_{k=0}^{n-1}\binom{n-1+k}{n-1-k}\frac{x^{k}}{k!}	\label{34}\\
&=n!\sum_{k=0}^{n}\binom{n+k-1}{n-k}\frac{x^{k}}{k!}. \nonumber
\end{align}
From \eqref{34}, we note that
\begin{align}
\sum_{n=0}^{\infty}L_{n}(x)\frac{t^{n}}{n!}&=\sum_{n=0}^{\infty}\bigg(\sum_{k=0}^{n}\binom{n+k-1}{n-k}\frac{x^{k}}{k!}\bigg)t^{n}=\sum_{k=0}^{\infty}\frac{x^{k}}{k!}\bigg(\sum_{n=k}^{\infty}\binom{n+k-1}{n-k}t^{n}\bigg)\label{35} 	\\
&=\sum_{k=0}^{\infty}\frac{x^{k}}{k!}\sum_{n=0}^{\infty}\binom{n+2k-1}{n}t^{n_k}=\sum_{k=0}^{\infty}\frac{(x)^{k}}{k!}\bigg(\frac{1}{1-t}\bigg)^{2k}=e^{\frac{xt}{(1-t)^{2}}}. \nonumber
\end{align}
Replacing $t$ by $1-e^{-t}$ in \eqref{35}, we have
\begin{align}
&\sum_{k=0}^{\infty}L_{k}(x)\frac{1}{k!}\big(1-e^{-t}\big)^{k}=e^{x(e^{2t}-1)}e^{-x(e^{t}-1)} \label{36}	\\
&=\sum_{k=0}^{\infty}\phi_{k}(x)2^{k}\frac{t^{k}}{k!}\sum_{m=0}^{\infty}\phi_{m}(-x)\frac{t^{m}}{m!}=\sum_{n=0}^{\infty}\bigg(\sum_{k=0}^{n}\binom{n}{k}\phi_{k}(x)2^{k}\phi_{n-k}(x)\bigg)\frac{t^{n}}{n!}. \nonumber
\end{align}
Proceeding in the same way as in \eqref{27}, we also have
\begin{equation}
\sum_{k=0}^{\infty}L_{k}(x)\frac{1}{k!}\big(1-e^{-t}\big)^{k}=\sum_{n=0}^{\infty}\bigg(\sum_{k=0}^{n}{n \brace k}(-1)^{n-k}L_{k}(x)\bigg)\frac{t^{n}}{n!}.\label{37}
\end{equation}
Therefore, by \eqref{36} and \eqref{37}, we obtain following theorem.
\begin{theorem}
For $n\ge 0$, we have
\begin{equation*}
\sum_{k=0}^{n}\binom{n}{k}2^{k}\phi_{k}(x)\phi_{n-k}(x)=\sum_{k=0}^{n}(-1)^{n-k}{n \brace k}L_{k}(x).
\end{equation*}
\end{theorem}
From \eqref{6} and \eqref{30}, we note that
\begin{align}
&\sum_{n=0}^{\infty}J_{n}(x)\frac{t^{n}}{n!}=\frac{1}{1-t}	e^{\frac{xt}{(1-t)^{2}}}=e^{-\log(1-t)}\sum_{k=0}^{\infty}L_{k}(x)\frac{t^{k}}{k!} \label{38} \\
&=\sum_{m=0}^{\infty}\frac{1}{m!}\log^{m}\bigg(\frac{1}{1-t}\bigg)\sum_{k=0}^{\infty}L_{k}(x)\frac{t^{k}}{k!}=\sum_{m=0}^{\infty}\sum_{j=m}^{\infty}{j \brack m}\frac{t^{j}}{j!}\sum_{k=0}^{\infty}L_{k}(x)\frac{t^{k}}{k!} \nonumber \\
&=\sum_{j=0}^{\infty}\sum_{m=0}^{j}{j \brack m}\frac{t^{j}}{j!}\sum_{k=0}^{\infty}L_{k}(x)\frac{t^{k}}{k!}=\sum_{n=0}^{\infty}\bigg(\sum_{j=0}^{n}\binom{n}{j}L_{n-j}(x)\sum_{m=0}^{j}{j \brack m} \bigg)\frac{t^{n}}{n!}. \nonumber
\end{align}
Therefore, by comparing the coefficients on both sides of \eqref{38}, we obtain the following theorem.
\begin{theorem}
For $n\ge 0$, we have
\begin{equation*}
J_{n}(x)=\sum_{j=0}^{n}\binom{n}{j}L_{n-j}(x)\sum_{m=0}
^{j}{j \brack m}.
\end{equation*}
\end{theorem}
From \eqref{34}, we determine $L_{n}(x)$, for the first few values of $n$:
\begin{align*}
&L_0(x) = 1, \quad L_1(x) = x, \quad L_2(x) = x^2 + 4x, \\
&L_3(x) = x^3 + 12x^2 + 18x, \quad L_4(x) = x^4 + 24x^3 + 120x^2 + 96x.
\end{align*}
We let $a_{0,n}(x)=L_{n}(x),\ (n\ge 0)$, in \eqref{2}. Then, by  the Seidel's formula in \eqref{12} and \eqref{35}, we get
\begin{equation*}
a_{n,0}(x)=\sum_{k=0}^{n}\binom{n}{k}L_{k}(x), \quad \sum_{n=0}^{\infty}a_{n,0}(x)\frac{t^{n}}{n!}=e^{t}e^{\frac{xt}{(1-t)^{2}}}.
\end{equation*}
The Euler-Seidel matrix associated with $(a_{0,n})_{n\ge 0}=\big(L_{n}(x)\big)_{n\ge 0}$ is given by
\begin{equation*}
\begin{aligned}
\left(\begin{matrix}
1 & x & x^{2}+4x & \cdots \\
x+1 & x^{2}+5x & x^{3}+13x^{2}+22x & \cdots \\
x^{2}+6x+1 & x^{3}+14x^{2}+27x & \vdots & \cdots \\
x^{3}+15x^{2}+33x+1 & \vdots & \vdots & \cdots\\
\vdots & \vdots & \vdots & \cdots
 \end{matrix}\right).
\end{aligned}
\end{equation*}

\section{Further Remark}
In this section, we study degenerate versions of $K_{n}(x)$, $J_{n}(x)$, and $L_{n}(x)$, namely $K_{n,\lambda}(x)$, $J_{n, \lambda}(x)$, and $L_{n,\lambda}(x)$, detailing their structural interrelations, their exponential generating  functions and their connection to the degenerate Euler–Seidel matrix framework. \par
First, we consider $K_{n,\lambda}(x)$ defined by
\begin{equation}
K_{n,\lambda}(x)=n!\sum_{k=0}^{n}\binom{n+k+1}{n-k}\frac{(x)_{k,\lambda}}{k!}, \quad (n \ge 0).\label{39}
\end{equation}
Then, by \eqref{4} and \eqref{39}, we get
\begin{align}
\sum_{n=0}^{\infty}K_{n,\lambda}(x)\frac{t^{n}}{n!}&=\sum_{n=0}^{\infty}\sum_{k=0}^{n}\binom{n+k+1}{n-k}\frac{(x)_{k,\lambda}}{k!}t^{n}\label{40}\\
&=\sum_{k=0}^{\infty}\frac{(x)_{k,\lambda}}{k!}\sum_{n=k}^{\infty}\binom{n+k+1}{n-k}t^{n}=\sum_{k=0}^{\infty}\frac{(x)_{k,\lambda}}{k!}t^{k}\sum_{n=0}^{\infty}\binom{n+2k+1}{n}t^{n} \nonumber\\
&=\sum_{k=0}^{\infty}\frac{(x)_{k,\lambda}}{k!}t^{k}\bigg(\frac{1}{1-t}\bigg)^{2k+2}=\frac{1}{(1-t)^{2}}e_{\lambda}^{x}\bigg(\frac{t}{(1-t)^{2}}\bigg). \nonumber
\end{align} \par
Second, we define $J_{n,\lambda}(x)$ by
\begin{equation}
J_{n,\lambda}(x)=K_{n,\lambda}(x)-nK_{n-1,\lambda}(x),\ (n\ge 0). \label{41}	
\end{equation}
Then, by \eqref{39} and \eqref{41}, we have
\begin{equation}
J_{n,\lambda}(x)=K_{n,\lambda}(x)-nK_{n-1,\lambda}(x)=n!\sum_{k=0}^{n}\binom{n+k}{n-k}\frac{(x)_{k,\lambda}}{k!}.\label{42}
\end{equation}
Thus, by \eqref{42}, we get
\begin{equation}
\sum_{n=0}^{\infty}J_{n,\lambda}(x)\frac{t^{n}}{n!}=\frac{1}{1-t}e_{\lambda}^{x}\bigg(\frac{t}{(1-t)^{2}}\bigg).\label{43}	
\end{equation}
From \eqref{40} and \eqref{43}, we note that
\begin{align}
\sum_{n=0}^{\infty}K_{n,\lambda}(x)\frac{t^{n}}{n!}&=\frac{1}{(1-t)^{2}}e_{\lambda}^{x}\bigg(\frac{t}{(1-t)^{2}}\bigg)=\frac{1}{1-t}\frac{1}{1-t}e_{\lambda}^{x}\bigg(\frac{t}{(1-t)^{2}}\bigg) \label{44} \\
&=\sum_{n=0}^{\infty}n!\sum_{k=0}^{n}\frac{J_{k,\lambda}(x)}{k!}\frac{t^{n}}{n!}.\nonumber
\end{align}
Therefore, by \eqref{44}, we obtain the following theorem.
\begin{theorem}
For $n\ge 0$, we have
\begin{equation*}
K_{n,\lambda}(x)=n!\sum_{k=0}^{n}\frac{1}{k!}J_{k,\lambda}(x).
\end{equation*}
\end{theorem}
Third, we define $L_{n,\lambda}(x)$ by (see \eqref{42})
\begin{equation}
L_{n,\lambda}(x)=J_{n,\lambda}(x)-nJ_{n-1,\lambda}(x)=n!\sum_{k=0}^{n}\binom{n-1+k}{n-k}\frac{(x)_{k,\lambda}}{k!},\quad (n \ge 0).\label{45}	
\end{equation}
Then, by \eqref{45}, we get
\begin{align}
\sum_{n=0}^{\infty}L_{n,\lambda}(x)\frac{t^{n}}{n!}&=\sum_{n=0}^{\infty}\sum_{k=0}^{n}\binom{n-1+k}{n-k}\frac{(x)_{k,\lambda}}{k!}t^{n} \label{46} \\
&=\sum_{k=0}^{\infty}\frac{(x)_{k,\lambda}}{k!}\sum_{n=k}^{\infty}\binom{n-1+k}{n-k}t^{n} =\sum_{k=0}^{\infty}\frac{(x)_{k,\lambda}}{k!}t^{k}\sum_{n=0}^{\infty}\binom{n+2k-1}{n}t^{n}\nonumber\\
&=\sum_{k=0}^{\infty}\frac{(x)_{k,\lambda}}{k!}\bigg(\frac{t}{(1-t)^{2}}\bigg)^{k}=e_{\lambda}^{x}\bigg(\frac{t}{(1-t)^{2}}\bigg). \nonumber
\end{align}
From \eqref{5}, \eqref{43} and \eqref{46}, we have
\begin{align}
&\sum_{n=0}^{\infty}J_{n,\lambda}(x)\frac{t^{n}}{n!}=\frac{1}{1-t}e_{\lambda}^{x}\bigg(\frac{t}{(1-t)^{2}}\bigg)=e_{\lambda}\bigg(\log_{\lambda}\bigg(\frac{1}{1-t}\bigg)\bigg)e_{\lambda}^{x}\bigg(\frac{t}{(1-t)^{2}}\bigg) \label{47} \\
&=e_{\lambda}\big(-\log_{-\lambda}(1-t)\big)e_{\lambda}^{x}\bigg(\frac{t}{(1-t)^{2}}\bigg)=\sum_{k=0}^{\infty}\frac{(1)_{k,\lambda}}{k!}(-1)^{k}\log_{-\lambda}^{k}(1-t)e_{\lambda}^{x}\bigg(\frac{t}{(1-t)^{2}}\bigg)\nonumber\\
&=\sum_{k=0}^{\infty}(1)_{k,\lambda}\sum_{m=k}^{\infty}{m \brack k}_{-\lambda}\frac{t^{m}}{m!}\sum_{j=0}^{\infty}L_{j,\lambda}(x)\frac{t^{j}}{j!}=\sum_{m=0}^{\infty}\bigg(\sum_{k=0}^{m}(1)_{k,\lambda}{m \brack k}_{-\lambda}\bigg)\frac{t^{m}}{m!}\sum_{j=0}^{\infty}L_{j,\lambda}(x)\frac{t^{j}}{j!} \nonumber\\
&=\sum_{n=0}^{\infty}\bigg(\sum_{m=0}^{n}\binom{n}{m}L_{n-m,\lambda}(x)\sum_{k=0}^{m}(1)_{k,\lambda}{m \brack k}_{-\lambda}\bigg)\frac{t^{n}}{n!}. \nonumber
\end{align}
Therefore, by comparing the coefficients on both sides of \eqref{47}, we obtain the following theorem.
\begin{theorem}
For $n\ge 0$, we have
\begin{equation*}
J_{n,\lambda}(x)=\sum_{m=0}^{n}\binom{n}{m}L_{n-m,\lambda}(x)\sum_{k=0}^{m}(1)_{k,\lambda}{m \brack k}_{-\lambda}.
\end{equation*}
\end{theorem}
From \eqref{45}, we find $L_{n,\lambda}(x)$, for the first few values of $n$:
\begin{align*}
&L_{0,\lambda}(x) = 1,\quad L_{1,\lambda}(x)= x, \quad L_{2,\lambda}(x) = x^2 + (4 - \lambda)x, \\
&L_{3,\lambda}(x) = x^3 + (12 - 3\lambda)x^2 + (2\lambda^2 - 12\lambda + 18)x, \\
&L_{4,\lambda}(x) = x^4 + (24 - 6\lambda)x^3 + (11\lambda^2 - 72\lambda + 120)x^2 \\
&\quad\quad\quad\quad+ (-6\lambda^3 + 48\lambda^2 - 120\lambda + 96)x.
\end{align*}
We let $a_{0,n}(x|\lambda)=L_{n,\lambda}(x),\ (n\ge 0)$, in \eqref{10}.
Then, by \eqref{13} and \eqref{46}, we get
\begin{align}
\sum_{n=0}^{\infty}a_{n,0}(x|\lambda)\frac{t^{n}}{n!}&=e_{\lambda}^{1-\lambda}(t)e_{\lambda}^{x}\bigg(\frac{t}{(1-t)^{2}}\bigg) \label{48}\\
&=\sum_{n=0}^{\infty}\sum_{k=0}^{n}\binom{n}{k}(1-\lambda)_{n-k,\lambda}L_{k,\lambda}(x)\frac{t^{n}}{n!}.\nonumber
\end{align}
Thus, by comparing the coefficients on both sides of \eqref{48}, we have
\begin{equation}
a_{n,0}(x|\lambda)=\sum_{k=0}^{n}\binom{n}{k}(1-\lambda)_{n-k,\lambda}L_{k,\lambda}(x). \label{49}
\end{equation}
The degenerate Euler-Seidel matrix associated with $\big(a_{0,n}(x|\lambda)\big)=\big(L_{n,\lambda}(x)\big)_{n\ge 0}$ is given by
\begin{align*}
\left(\begin{matrix}
1 & x  & \cdots \\
x+(1-\lambda) & x^{2}+(5-\lambda)x  & \cdots \\
x^{2}+(6-3\lambda)x+(1-\lambda)_{2,\lambda} & x^{3}+(14-3\lambda)x^{2}+(27-15\lambda+2\lambda^{2})x & \cdots \\
\vdots & \vdots & \vdots
 \end{matrix}\right).
\end{align*}
\section{Conclusion}
In this study, we successfully combined the algebraic utility of the Euler–Seidel matrix method with the structural framework of degenerate polynomials and Morgan–Voyce type polynomials. \\
$\indent \bullet$ New Polynomial Families: We introduced and systematically characterized several new Morgan–Voyce type polynomials, namely $M_{n}(x)$, $N_{n}(x)$, $K_{n}(x)$, $J_{n}(x)$, and $L_{n}(x)$, deriving their explicit expressions, exponential generating functions and recurrence relations. \\
$\indent \bullet$ Binomial-type Identities: We established new connections to Bell polynomials, proving that three separate binomial-type sums are uniquely equal to finite sums of our newly defined variants. \\
$\indent \bullet$ Degenerate Extensions: Moving into the degenerate domain, we constructed the parameter-dependent variants $K_{n,\lambda}(x)$, $J_{n,\lambda}(x)$, and $L_{n,\lambda}(x)$, mapping out their corresponding properties and interactions via the degenerate Euler–Seidel matrix.\par
The identities and relations uncovered in this work demonstrate that integrating matrix transformation methods with degenerate polynomial frameworks serves as a fertile ground for discovering hidden structures in combinatorial analysis and analytic number theory. Future research may extend these methodologies to other classical special polynomials, or explore the $p$-adic properties of the degenerate polynomial families introduced herein.

\vspace{0.5cm}

 {\bf Competing interests.} The authors have no competing interests to declare that are relevant to the content of this article and the authors contributed equally to this work.

\end{document}